\documentclass [12pt]{article}

\usepackage{amsmath, amssymb,amsgen, amsthm, amscd}

\newtheorem{thm}{Theorem}
\newtheorem{lemma}[thm]{Lemma}
\newtheorem{cor}[thm]{Corollary}
\newtheorem{prop}[thm]{Proposition}

\theoremstyle{definition}

\theoremstyle{remark}
\newtheorem{rem}{Remark}

\numberwithin{equation}{section}
\numberwithin{thm}{section}
\numberwithin{rem}{section}

\newcommand{\ca}{\mathcal{A}}

\newcommand{\cc}{\mathcal{C}}
\newcommand{\cg}{\mathcal{G}}
\newcommand{\cx}{\mathcal{X}}
\newcommand{\ch}{\mathcal{H}}

\newcommand{\lp}{\mathcal{L}^p}

\newcommand{\tcc}{\widetilde{\mathcal{C}}}
\newcommand{\tca}{\widetilde{\mathcal{A}}}
\newcommand{\tWW}{\widetilde{W}}
\newcommand{\tf}{\widetilde{f}}

\newcommand{\al}{\alpha}
\newcommand{\f}{\phi}

\newcommand{\de}{\delta}
\newcommand{\De}{\Delta}

\newcommand{\w}{\omega}

\newcommand{\W}{\Omega}
\newcommand{\tW}{\widetilde{\Omega}}
\newcommand{\bt}{\beta}
\newcommand{\g}{\gamma}
\newcommand{\G}{\Gamma}

\newcommand{\s}{\sigma}
\newcommand{\na}{\nabla}
\newcommand{\tna}{\widetilde{\nabla}}

\newcommand{\te}{\theta}
\newcommand{\Te}{\Theta}
\newcommand{\tte}{\widetilde{\theta}}

\newcommand{\tch}{\widetilde{\ch}}
\newcommand{\tF}{\widetilde{F}}
\newcommand{\tg}{\widetilde{\g}}
\newcommand{\del}{\partial}
\newcommand{\df}{\frac{dF_t}{dt}}
\newcommand{\gotimes}{\widehat{\otimes}}
\newcommand{\dint}{{\displaystyle \int \!\!\!\!\!\!-}}

\newcommand{\cp}{\sqrt{2\pi i}}
\newcommand{\ci}{\sqrt{2 i}}
\newcommand{\ccc}{\mathcal{C}^c}
\newcommand{\Wc}{\W^c}
\newcommand{\nac}{\na^c}
\newcommand{\tec}{\te^c}
\newcommand{\intc}{\dint^c}

\DeclareMathOperator{\ind}{index}
\DeclareMathOperator{\ad}{ad}
\DeclareMathOperator{\End}{End}
\DeclareMathOperator{\codim}{codim}
\DeclareMathOperator{\Diff}{Diff}
\DeclareMathOperator{\Ch}{Ch}
\DeclareMathOperator{\flow}{sf}

\begin{document}
\title{Characters of  Cycles, Equivariant Characteristic Classes
 and Fredholm Modules}
\author{Alexander Gorokhovsky
\\Department of Mathematics,\\
 Ohio State
University,\\
Columbus, OH 43210\\
sasha@math.ohio-state.edu}
\maketitle

\begin{abstract}
We derive simple explicit formula for the character of a cycle in the
 Connes' $(b,$ $B)$-bicomplex of cyclic cohomology
 and apply it to write formulas for the equivariant Chern character and
characters of finitely-summable bounded Fredholm modules.
\end{abstract}

\section{Introduction}
 The  notion of a cycle,  introduced by Connes in \cite{cn85},
 plays an important role in his  development of
the cyclic cohomology  and its applications.   Many
questions of the  differential geometry and noncommutative geometry
can be reformulated as questions about  geometrically  defined cycles.
Associated with a cycle is its \emph{character}, which is a
characteristic class in cyclic cohomology,
described by an explicit formula ( see \cite{cn85}).

%The definition of a cycle imposes however a certain ``flatness'' condition.
% If one removes it one arrives at the definition
%of what we call a ``generalized cycle''.
% Such objects appear naturally in Connes' work
%(\cite{cn85},\cite{cn94}). In the context of superconnections they
%were studied in \cite{ni95}.
%To deal with them Connes has devised a canonical procedure for annihilating the curvature. 
Some natural constructs, like the
 the transverse fundamental cycle of a foliation \cite{cn94} or
the superconnection  in \cite{ni95} require
however  consideration of  more general objects, which we call
``generalized cycles'' (we recall the definition in the section \ref{cycles}).
The simplest geometric example of generalized cycle is provided by the
algebra of forms with values in the endomorphisms of some vector bundle,
together with a connection. More interesting examples arise from 
vector bundles equivariant with respect to action of discreet group, or,
more generally, holonomy equivariant vector bundles on foliated manifolds.

The original definition of the character of a cycle does not apply  directly 
to generalized cycles.
To overcome this,    Connes (\cite{cn85}, cf. also \cite{cn94})
has devised a canonical procedure
allowing to associate a cycle with a generalized cycle. This allows to
extend the definition of the character to the generalized cycles.

In this paper we  show that  the character of a generalized cycle
can be defined by the 
explicit  formula  in the
$(b,$ $B)$-bicomplex, resembling the JLO formula for the Chern character
\cite{jlo88}. In the geometric examples above this leads to
formulas  for the Bott's Chern character \cite{bott78} in cyclic
cohomology. As another example we derive formula for the character of Fredholm
module.

The paper is organized as follows.
In   section \ref{cycles} we define the character of a generalized cycle and,
more generally, generalized chain.
 Closely related formulas 
also appear and play an important role in Nest and Tsygan's work on the
algebraic index theorems \cite{nt95, nt95a}.
We then establish some basic properties of this character and prove that our
definition of the character coincides with the original
one given by Connes
in \cite{cn85}.
In   section \ref{equiv} we 
construct the cyclic cocycle,
representing the equivariant  Chern character
in the cyclic cohomology,
and discuss relation of this construction
 with the multidimensional version of the Connes
construction of the Godbillon-Vey cocycle \cite{cn86}, and the
transverse fundamental class of the foliation.
In   section \ref{appl} we
  write explicit formulas for the character of a
bounded finitely-summable
Fredholm module, where $F^2-1$ is not necessarily 0~
 (such objects are called
pre-Fredholm modules in \cite{cn85}). The idea is to associate with
 such a Fredholm module 
 a generalized cycle, by the construction similar to
\cite{cn85}. We thus obtain finitely summable analogues of the  formulas from
\cite{jlo88} and \cite{gs89}.

I would like to thank my advisor H. Moscovici for introducing me to the
area and constant support.
I would like also to thank D. Burghelea and I. Zakharevich
for helpful discussions.

\section{Characters of cycles} \label{cycles}
In this section we start by stating definitions of generalized chains and
cycles, and  writing the JLO-type formula for the character. We then
show that this definition of character coincides with the original one
from \cite{cn85}.

In what follows we require the algebra $\ca$ to be unital. This
condition will be later removed by adjoining the unit to $\ca$.

%A generalized cycle over an algebra $\ca$  is given by the following
%data :
%\begin{enumerate}
%\item A $\mathbb{Z}$-graded unital algebra
%$\W=\bigoplus_{m=0}^{\infty}\W^m$ and a homomorphism $\rho$ from $\ca$ to
%$\W^{\,0}$. We require the homomorphism to be unital.

%\item A graded derivation $\na:\W^k\mapsto \W^{k+1}$, $k=$0~, 1~,\dots
% and $\te\in \W^2$ such that            
% \[\na^2(\xi)=\te\xi-\xi\te\] $\forall \xi \in \W$. We require that
%$\na(\te)=0$ . We will call $\na$ connection and $\te$
%curvature .

%\item \label {trace}A graded trace $\dint$ defined on $\W^{n}$ for some  $n$ with the
%property \[\dint \na (\xi)=0 \] $\forall \xi \in \W^{n-1}$. We will call
%$n$ a degree of a cycle.
%\end{enumerate}
%Note that for $n$ odd notions of a graded trace and a trace coincide. 
%The definition of cycle is obtained by requiring $\te$ to be 0~.

One  defines a \emph{generalized chain} over an algebra $\ca$
by specifying the following
data:
\begin{enumerate}
\item Graded unital algebras $\W$ and $\del \W$ and a surjective homomorphism $r:
\W \rightarrow \del \W$ of degree 0~, and a homomorphism $\rho :\ca
\rightarrow \W^0$. We require that $\rho$ and $r$ be unital.

\item Graded derivations of degree 1~ $\na$ on $\W$ and $\na'$ on $\del
\W$ such that $r \circ \na=\na' \circ r$ and  $\te \in
\W^2$ such that 
\[
\na^2(\xi)=\te\xi-\xi\te
\]
$\forall \xi \in \W$. We require that $\na(\te)=0$ .

\item A graded trace $\dint$ on $\W^n$ for some $n$ (called the degree of
the chain) such that
\[
\dint \na (\xi)=0
\]
$\forall \xi \in \W^{n-1}$ such that $r(\xi)=0$.

\end{enumerate}

If one requires $\del \W=0$ one obtains the definition of the
\emph{generalized cycle}. Generalized cycle for which $\te=0$ is called
\emph{cycle}. 

One defines the boundary of the generalized chain to be a generalized
cycle $(\del \W, \na', \te',\dint')$ of degree $n-1$ over an algebra $\ca$
where the $\dint'$ is the graded trace
defined by the identity
\begin{equation}\label{defint'}
\dint' \xi'=\dint \na(\xi)
\end{equation} where $\xi'\in (\del
\W)^{n-1}$ and $\xi \in \W^n$ such that $r(\xi)=\xi'$. Homomorphism
$\rho': \ca \rightarrow \del\W^0$ is given by
\begin{equation} \label {defr'}
\rho'=r\circ \rho
\end{equation}
Notice that for $\xi'\in \del \W$ $(\na')^2(\xi')=\te' \xi'-\xi' \te'$
where $\te'$ is defined by
\begin{equation} \label{defte'}
\te'=r(\te)
\end{equation}

With every generalized chain $\cc^n$ of degree $n$
one can associate by a JLO-type formula a canonical
$n$-cochain $\Ch(\cc^n)$ in the $(b,$ $B)$-bicomplex of the algebra $\ca$,
which we call a character of the generalized chain.

\begin{multline}\label{dchern}
\Ch^k(\cc^{n})(a_0,a_1,\dots a_k)=\\
\frac{(-1)^{\frac{n-k}{2}}}{(\frac{n+k}{2})!}\sum_{i_0+i_1+
\dots + i_k=\frac{n-k}{2}}\dint \rho(a_0) \te^{i_0} \na(
\rho(a_1))\te^{i_1}\dots \na(\rho(a_k)) \te^{i_k}
\end{multline}

Note that if $\cc^n$ is a (non-generalized) cycle $\Ch(\cc^n)$ coincides
with the character of $\cc^n$ as defined by Connes. 

For the generalized chain $\cc$ let $\del \cc$ denote the boundary of $\cc$.
\begin{thm} \label{bch}
Let $\cc^n$ be a chain, and $\del(\cc^n)$ be its boundary. Then
\begin{equation}
(B+b)
\Ch(\cc^{n})=     
  S\,\Ch(\del(\cc^{n})) 
\end{equation}
Here $S$ is the usual periodicity shift in the cyclic bicomplex.
\end{thm}
\begin{proof} By direct computation.
\end{proof}
\begin{rem}
A natural framework for such identities in cyclic cohomology is
provided by the theory of operations on cyclic cohomology of Nest and Tsygan,
cf. \cite{nt95,nt95a} 
\end{rem}
\begin{cor} \label{cocycle}
If $\cc^n$ is a generalized cycle then $\Ch(\cc^n)$ is an $n$-cocycle in the
cyclic bicomplex of an algebra $\ca$. 
\end{cor}
\begin{cor} \label{cobord}
 For two cobordant  generalized cycles $\cc^n_1$ and $\cc^n_2$   
\[
[S\,\Ch(\cc^n_1)]=[S\,\Ch(\cc^n_2)]
\]
 in $HC^{n+2}(\ca)$.
\end{cor}

 Formula \eqref{dchern} can also be written in the
different form. We will use the following notations.
First, $\dint$ can be
extended to the whole algebra $\W$ by setting $\dint \xi=0$ if $\deg \xi
\ne n$.
For $\xi \in \W$ $e^{\xi}$ is defined as $\sum_{j=0}^{\infty} \frac{\xi^j}{j!}$ .
Then denote $\De^k$ the $k$-simplex $\{(t_0, t_1, \dots,
t_k)|t_0+t_1+\dots+t_k=1, t_j \geq 0 \}$ with the measure $dt_1dt_2\dots dt_k$.
Finally, $\al$ is an arbitrary nonzero real parameter.
Then
\begin{multline}\label{ichern}
\Ch^k(\cc^{n})(a_0,a_1,\dots a_k)=\\
 \al^{\frac{k-n}{2}} \int \limits_{\De^k}\left( \dint \rho(a_0) e^{-\al t_0\te} \na(
\rho(a_1)) e^{-\al t_1\te} \dots \na(\rho(a_k)) e^{-\al t_k\te}  \right)dt_1dt_2\dots dt_k
\end{multline}
where
 $k$ is of the same parity as $n$.
Indeed,
\begin{multline}
\dint \rho(a_0) e^{-\al t_0\te} \na(
\rho(a_1)) e^{-\al t_1\te} \dots \na(\rho(a_k)) e^{-\al t_k\te}=\\
(-\al)^{\frac{n-k}{2}}\sum_{i_0+i_1+
\dots + i_k=\frac{n-k}{2}} \frac{t_0^{i_0}t_1^{i_1}\dots t_k^{i_k}}{i_0!i_1!
\dots i_k!}\dint \rho(a_0) \te^{i_0} \na(
\rho(a_1))\te^{i_1}\dots \na(\rho(a_k)) \te^{i_k}
\end{multline}
and our assertion follows from the equality
\[
\int \limits_{\De^n} t_0^{i_0}t_1^{i_1}\dots t_k^{i_k} \,dt_1dt_2\dots dt_k=
\frac{i_0!i_1!\dots i_k!}{(i_0+i_1+\dots +i_k+k)!}
\]

\begin{rem}\label{unit}

We worked above only in the context of unital algebras and maps.
The
case of general algebras and maps can be treated by adjoining a unit.
 We follow \cite{ni95}
The definition of the generalized chain in the nonunital case differ
from the definition in the unital case only in two aspects:
first, we do not require algebras and morphisms to be unital,
second, we do not require any more that the curvature $\te$ is an
element of $\W^2$; rather we require it to be a multiplier of the
algebra $\W$ which satisfies  the following:
for $\w \in \W^k$ $\te \w$ and $\w \te$ are in $\W^{k+2}$, 
$\na(\te \w)=\te \na(\w)$, $\na(\w \te) =\na(\w) \te$ and  $\dint \te \w=\dint
\w \te$ if $\w \in \W^{n-2}$. We also need to require existence of
$\te'$ -- multiplier of
$\del \W$ such that $r(\te \w)=\te' r(\w)$, $r(\w \te)= r(\w)\te'$, and
include it in the defining data of chain.

With $\cc^n=(\W,\del \W,r,\na,\na', \te, \dint)$ -- nonunital generalized
chain over a (possibly nonunital) algebra $\ca$
we  
associate canonically  a unital chain
$\tcc^n=(\tW,\del \tW,\widetilde{r},\tna,\tna'\tte, \widetilde{\dint})$
over the algebra $\tca$ -- $\ca$ with unit adjoined.
The construction is the following:
the algebra $\tW$ 
is obtained from the algebra $\W$ 
by adjoining a unit 1~,
( of degree 0~) and an element $\tte$ of degree 2~ with the
relations $\tte \w=\te \w$ and $\w \te= \w \tte$ for $\w \in \W$,
and similarly for the algebra $\del \tW$.
The derivation $\tna$ coincides with $\na$ on the elements of
$\W$ and satisfies equalities $\tna(\tte)=0$ and $\tna(1)=0$,
and  $\tna'$ is defined similarly.
The graded
trace $\widetilde{\dint}$ on $\tW$ is defined to coincide with $\dint$ on the elements
of $\W$ and, if $n$ is even, is required to satisfy the relation $\widetilde{\dint}
\tte^{\frac{n}{2}} =0$.

Now if $\cc^n$ is a (nonunital) generalized cycle over
 $\ca$, formula \eqref{dchern}, applied to to $\tcc^n$  defines a
(reduced) cyclic cocycle over an algebra $\tca$ and hence
 a class in the reduced cyclic cohomology
$\overline{HC}^n(\tca) =HC^n(\ca)$. The Corollary \ref{cobord}
implies that this class is invariant under the (nonunital) cobordism.
Note also that in the unital case the class defined after adjoining the
unit agrees with the one defined before.

Alternatively, one can work from the beginning
with the Loday-Quillen-Tsygan bicomplex,
see e.g. \cite{lo92}, where the corresponding formulas can be easily written.

\end{rem}

We now will show equivalence of the previous construction with Connes'
original construction.

With every generalized cycle $\cc=(\W,\na,\te,\dint)$ over an algebra $\ca$
Connes shows how to associate canonically a
cycle $\cc_X$.

One starts with a
 graded algebra $\W_{\te}$, which as a vector space can be identified with
the space of 2~ by 2~ matrices over an algebra $\W$, with the grading
given by the following:
\[
 \left[
 \begin{matrix}
 \w_{11}\ \w_{12} \\
 \w_{21}\ \w_{22}
 \end{matrix}
 \right] \in  \W_{\te}^k
\text{ if }
\w_{11}\in \W^k \ 
\w_{12}, \w_{21}\in \W^{k-1} \text{ and }
\w_{22}\in \W^{k-2}
\]

The product of the two elements in $\W_{\te}$ $\w=\left[
 \begin{matrix}
 \w_{11}\ \w_{12} \\
 \w_{21}\ \w_{22}
 \end{matrix}
 \right]
$
and $\w'=\left[
 \begin{matrix}
 \w'_{11}\ \w'_{12} \\
 \w'_{21}\ \w'_{22}
 \end{matrix}
 \right]
$
is given by
\begin{equation}
\w *\w'=\left[
 \begin{matrix}
 \w_{11}\ \w_{12} \\
 \w_{21}\ \w_{22}
 \end{matrix}
 \right]
 \left[
 \begin{matrix}
 1 \quad 0 \\
 0 \quad  \te
 \end{matrix}
 \right]
\left[
 \begin{matrix}
 \w'_{11}\ \w'_{12} \\
 \w'_{21}\ \w'_{22}
 \end{matrix}
 \right]
\end{equation}
The homomorphism $\rho_{\te}:\ca \rightarrow \W_{\te}$ is given by
\begin{equation}
\rho_{\te}(a)=
\left[
 \begin{matrix}
 \rho(a) &0\\
 0       &0
 \end{matrix}
\right]
\end{equation}
On this algebra one can define a graded derivation $\na_{\te}$ of degree
1~ by the formula (here $\w=\left[
 \begin{matrix}
 \w_{11}\ \w_{12} \\
 \w_{21}\ \w_{22}
 \end{matrix}
 \right]
$)
\begin{equation}
\na_{\te} (\w)=
 \left[
 \begin{matrix}
 &\na (\w_{11})\ &\na(\w_{12}) \\
 &-\na(\w_{21})\ &-\na(\w_{22})
 \end{matrix}
 \right]
\end{equation}
One checks that
\begin{equation}\label{curv}
\na_{\te}^2 (\w)=
\left[
\begin{matrix}
\te \quad 0 \\
0 \quad 1
\end{matrix}
\right]*\w-
\w*\left[
\begin{matrix}
\te \quad 0 \\
0 \quad 1
\end{matrix}
\right]
\end{equation}
More generally, one can define on this algebra a family of
connections $\na^t_{\te}$, $0\leq t \leq 1$ by
the equation
\begin{equation}
\na^t_{\te}(\w)=\na_{\te}(\w) + t(\cx*\w -(-1)^{\deg\w}\w* \cx)
\end{equation}
where $\cx$ is degree 1~ element of $\W_{\te}$
given by the matrix
\begin{equation}
\cx=
\left[
\begin{matrix}
0 \ &-1\\
1 \ &0
\end{matrix}
\right]
\end{equation}
\begin{lemma}
$(\na_{\te}^t)^2(\w)=(1-t^2)
\left(
\left[
\begin{matrix}
 \te \quad 0\\
 0 \quad 1
 \end{matrix}
\right]*\w -\w*\left[
\begin{matrix}
 \te \quad 0\\
 0 \quad 1
 \end{matrix}
\right]
\right)
$
\end{lemma}
\begin{proof}
%We have: $\na_{\te}^t=\na_{\te}+t \ad_{\cx}$ (by $\ad_{\cx}$
%we mean the action by the graded commutator). Since $\na_{\te}(\cx)=0$,
%$\na_{\te}$ and $ \ad_{\cx}$ anticommute, and 
%\[
%(\na_{\te}^t)^2=(\na_{\te})^2+t^2 (\ad_{\cx})^2=(\na_{\te})^2+t^2 \ad_{\cx^2}
%\]
%and the assertion follows from the \eqref{curv} and the identity
%\[
%\cx^2=\cx * \cx=
%\left[
%\begin{matrix}
% &-\te \quad &0\\
% &0 \quad &-1
% \end{matrix}
%\right]
%\]
Follows from an easy computation.
\end{proof}
 Hence for $t=1$ we obtain a graded derivation $\na_{\te}^1$ whose
square is 0~.

Finally, the graded trace $\dint_{\te}$ is defined by
\begin{equation}
\dint_{\te}\w=\dint \w_{11}-(-1)^{\deg \w}\dint \w_{22}\te
\end{equation}
It is closed with respect to $\na_{\te}$, and hence, being a graded
trace, it is closed with respect to $\na_{\te}^t$ for any $t$.

\begin{cor}
$\cc_X=(\W_{\te}, \na_{\te}^1, \dint_{\te})$ is a (nongeneralized) cycle
\end{cor}
The cycle $\cc_X$ is the Connes' canonical cycle, associated with the
generalized cycle $\cc$ .
With every (nongeneralized) cycle  of degree $n$ Connes associated a
cyclic $n$-cocycle on the algebra $\ca$ by
the following procedure: let the cycle consist of a graded
 algebra $\W$, degree 1~ graded derivation $d$ and a closed trace
$\dint$. Then the character of the cycle is the cyclic cocycle $\tau$ in the
 \emph{cyclic complex} 
given by the formula
\begin{equation}
\tau( a_0,a_1,\dots, a_n)=\dint \rho(a_0) d\rho(a_1) \dots d \rho(a_n)
\end{equation}
To it corresponds a cocycle in the $(b,$ $B)$-bicomplex with only the one
nonzero component of degree $n$, which equals
$
\frac{1}{n!}\dint \rho(a_0) d\rho(a_1) \dots d \rho(a_n)
$
\begin{thm}\label{comp}
Let $\cc^n$ be a generalized cycle of degree $n$ over an algebra $\ca$, 
 and $\cc_X^n$ be the canonical cycle over $\ca$, associated with $\cc^n$
 (see above). Then
$[\Ch(\cc^n)]=[\tau(\cc_X^n)]$ in $HC^n(\ca)$.
\end{thm}
Note that equality here is in the cyclic cohomology, not only in the
periodic cyclic cohomology.
The theorem will follow easily from the above considerations and the
following lemma. 
\begin{lemma}\label{varcon}
Let $(\W, \na_t, \te_t, \dint)$ ,$0\leq t \leq 1$ be a family of cycles
of degree $n$
over an algebra $\ca$
with  connection and curvature depending 
on $t$, all the other data
staying the same. Connection and curvature vary by
\[
\na_t= \na_0 +t \ad \eta
\]
for
some $\eta \in \W^1$, and
\[
\te_t=\te_0 +t \na_0 \eta +t^2 \eta^2
\]
Let $\Ch(t)$ denote the cocycle  obtained for some
specific value of $t$. Then $[\Ch(0)]=[\Ch(1)]$.
\end{lemma}
\begin{proof}[Proof of the Lemma \ref{varcon}]
First, we can suppose that the cycle is unital -- in the other case one
can perform a construction, explained in the Remark \ref{unit}.
 We define $\widetilde{\te_t}$ by
$  \widetilde{\te_t}=\widetilde {\te_0} +t \na_0 \eta +t^2 \eta^2$. 

We start by constructing a cobordism between cycles obtained when $t=0$
and $t=1$. This is analogous to a construction from \cite{ni95}.
The cobordism is provided by the chain $\ccc$ defined as follows:
The algebra $\Wc=\W^*([0,1]) \gotimes \W$, where $\gotimes$ denotes the
graded tensor product, and $\W^*([0,1])$ is the algebra of the
differential forms on the segment $[0,1]$.
The map $\rho^c:\ca \rightarrow \Wc$ is given by
\begin{equation}
\rho^c(a)=1 \gotimes \rho(a)
\end{equation}
We denote by $t$ the
variable on the segment $[0,1]$.
On this algebra we introduce
the graded derivation
\begin{equation}
\nac=1 \gotimes \na_t +d \gotimes 1
\end{equation}
 Here $d$ is
the de Rham differential, and $1 \gotimes \na_t$ has the following
meaning: for $\w \in \W$ $\na_t(\w)$ can be considered as an element in
$C^{\infty}([0,1]) \otimes \W^*([0,1]) \subset \W\gotimes \W$. We then
define (for $\w\in \W$, $\al \in \W^*([0,1])$ )
\[
1 \gotimes \na_t \,(\al\gotimes \w)=(-1)^{\deg \al} (\al\gotimes 1)\na_t(\w)
\]
The curvature $\tec$ is defined to be $1\gotimes \te_t+dt \gotimes
\eta$.  The restriction map $r^c$ maps $\Wc$ to $\W \oplus \W$, and
given by the restriction of function to the interval endpoints. Here we
consider the first $\W$ to come from the cycle obtained for $t=1$ and
the second from the cycle, obtained for $t=0$ and the trace given by
$-\dint$.
 
The graded trace $\intc$ on $(\Wc)^{n+1}$ is given by the formula
\begin{equation}
\intc(\al \gotimes \w)=
\begin{cases}
 (-1)^n\dint\w \int \limits_{[0,1]}\al &\text{ if } \deg \w=n \text { and
} \deg\al =1\\
0  &\text{ otherwise}
\end{cases}
\end{equation}
One checks that the above construction gives a chain providing cobordism
between cycles obtained when $t=0$ and $t=1$. Notice also that
the Theorem \ref{bch} provides an explicit cochain $\Ch(\ccc)$ of degree
$n+1$ such that $(b+B)\, \Ch(\ccc)=S\,\left(\Ch(1)-\Ch(0)\right)$. Its top
component, $\Ch^{n+1}(\ccc)$ is given by the formula
\[
\intc \rho^c(a_0) \nac(\rho^c(a_1))\dots \nac(\rho^c(a_{n+1}))
\]
but this  is easily seen to be 0~, since the expression under
the $\intc$ does not contain $dt$. This means that we can consider
$\Ch(\ccc)$ as a $n-1$ chain, for which then $(b+B)\, \Ch(\ccc)=\Ch(1)-\Ch(0)$,
and this proves the Lemma.
\end{proof}
\begin{rem} The above lemma remains true if we relax its conditions to allow
$\eta$ to be a multiplier of degree 1~, such that $\dint \eta \w =
(-1)^{(n-1)/2}\dint \w \eta$. Then $ \na_0 \eta$ is a multiplier,
 defined by $(\na_0 \eta) \w = \na_0 (\eta\w)
+ \eta \na_0\w $. The same proof then goes through if we enlarge the
algebra $\W$ to the subalgebra of the multiplier algebra of $\W$
obtained from $\W$ by adjoining 1~, $\te_0$, $\eta$, $\na_0 \eta$, and
extending $\dint$ to this algebra by zero (i.e. we put $\dint P =0$ for
any $P$ --  monomial in $\te_0$ and $\eta$).
\end{rem}

\begin{proof}[Proof of the Theorem \ref{comp}]
 The lemma above applies to the family
of cycles constructed above (with the curvatures 
$\te_t=
(1-t^2)
\left[
\begin{matrix}
\te& 0\\
0& 1
\end{matrix}
\right]
$ and connections $\na_t=\na_{\te}^t$).
This implies
 that $\Ch(\cc^n)=\Ch(\cc_X^n)$ in $HC^n(\ca)$.
Since $\cc_X^n$ is a (nongeneralized) cycle, comparison of the
definitions shows that  $\Ch(\cc_X^n)= \Ch(\cc_X^n)$, even on
the level of cocycles, and the Theorem follows.
\end{proof}

\begin{cor} \label{prod}
For two generalized cycles $\cc_1^n=(\W_1, \na_1,\te_1, \dint_1)$ and
$\cc_2^m=(\W_2, \na_2,\te_2,
\dint_2)$ define the product by $\cc_1\times \cc_2 =(\W_1\gotimes \W_2, \na_1
\gotimes 1 +1\gotimes \na_2,\te_1\gotimes 1+1 \gotimes \te_2 , \dint_1 \gotimes \dint_2)$.
Then $\Ch(\cc_1\times \cc_2)= 
\Ch(\cc_1)\cup \Ch(\cc_2)$.
\end{cor}
\begin{proof}
For the non-generalized cycles this follows
from  Connes' definition of the
cup-product. In the general case,
the statement follows from the existence of the natural map of cycles
(i.e. homomorphism of the corresponding algebras, preserving all the structure)
$(\cc_1\times \cc_2)_X \to (\cc_1)_X\times (\cc_2)_X$, which agrees with
taking of the character.

The simplest way to describe this map is by
using another Connes' description of his construction. In this
description matrix $\left[ \begin{matrix}&\w_{11} \quad &\w_{12}\\
&\w_{21} \quad &\w_{22} \end{matrix} \right]$, $\w_{ij} \in \W$ is
identified with the element $\w_{11}+ \w_{12}X +X\w_{21} + X\w_{22}X$,
where $X$ is a formal symbol of degree 1~. The multiplication law  is
formally defined by $\w X \w'=0$, $X^2=\te$. This should be understood
as a short way of writing identities like $\w X * X\w' =\w\te \w'$ (note
that $X$ is not an element of the algebra).

If we denote by $X_1$, $X_2$, $X_{12}$ formal elements, corresponding to $\cc_1$,
$\cc_2$, $\cc_1\times \cc_2$ respectively, the homomorphism 
mentioned above is
the unital
extension  of the identity map $\W_1\gotimes \W_2 \to \W_1\gotimes
\W_2$ defined ( again formally) by $X_{12}\mapsto (X_1\gotimes1+1\gotimes X_2)$.

\end{proof}

\section{Equivariant characteristic classes}\label{equiv}
This section concerns vector bundles equivariant with respect to
discrete group actions. We show that there is a generalized cycle
associated naturally to such a bundle with ( not necessarily invariant )
connection. The character of this generalized cycle turns out to be
related ( see Theorem \ref{phi} ) to the equivariant Chern character.

 Let $V$ be an orientable
 smooth manifold of dimension $n$, $E$ a complex vector bundle over $V$,
and $\ca=\End(E)$ -- algebra of endomorphisms
with compact support. One can construct a generalized cycle over an algebra
$\ca$ in the following way. The algebra $\W=\W^*(V,\End(E))$ -- the
algebra of endomorphism-valued differential forms. Any connection $\na$
on the bundle $E$ defines a connection for the generalized cycle, with
the curvature $\te \in \W^2(V,\End(E))$ -- the usual curvature of the
connection. On the $ \W^n(V,\End(E))$ one defines a graded trace $\dint$ by
the formula $\dint \w=\int \limits_V tr \w$, where in the right hand side
we have a usual matrix trace and a usual integration over a manifold.
Note that when $V$ is noncompact, this cycle is nonunital. 
The formula \eqref{ichern},  define a cyclic
$n$-cocycle $\{\Ch^k\}$ on the algebra $\ca$, given by the formula
\begin{multline}
\Ch^k(a_0,a_1,\dots a_k)=\\
  \int \limits_{\De^k}\left( \int \limits_V tr\, a_0  e^{- t_0\te} \na(
a_1) e^{-t_1\te} \dots \na(a_k) e^{-t_k\te}  \right)dt_1dt_2\dots dt_k
\end{multline}
Hence we recover the formula of Quillen from \cite{qu88}.
(Recall that for noncompact $V$ these expressions should be viewed as
defining reduced cocycle over the algebra $\ca$ with unit adjoined, with
$\Ch^0$ extended by $\Ch^0(1)=0$).

One can restrict this cocycle to the subalgebra of functions
$C^{\infty}(V) \subset \End(E)$. As a result one obtains an $n$-cocycle on the
algebra $C^{\infty}(V)$, which we still denote by $\{\Ch^k\}$, given by the formula
\begin{equation}
\Ch^k(a_0,a_1,\dots a_k)=
  \frac{1}{k!} \int \limits_V  a_0   
da_1  \dots da_k \, tr\, e^{-\te}  
\end{equation}
To this cocycle corresponds a current on $V$, defined by the form
$ tr\,e^{-\te}$.  Hence in this case
we recover the Chern character of the bundle $E$. Note
that we use normaliztion of the Chern character from \cite{bgv}.

 Let now an orientable manifold $V$ of dimension $n$ be equipped with
 an action of the
discreet group $\G$ of orientation preserving transformations, and $E$
be a $\G$-invariant bundle. In this situation, one can again construct a
cycle of degree $n$ over the algebra $\ca=\End(E) \rtimes \G$.
Our notations are the following : the algebra $\ca$ is generated by the
elements of the form $a U_g$, $a\in \End(E)$, $g\in \G$, and $U_g$ is a
formal symbol.
The
product is $(a'U_
{g'})(a U_g)=a'a^{g'} U_{gg'}$. The
superscript here denotes the action of the group.

The graded
algebra $\W$ is defined as $\W^*(V,\End(E))\rtimes \G$.
Elements
of $\W$ clearly act  on the forms with values in $E$,
and any connection $\na$ in the
bundle $E$ defines a connection for the algebra $\W$, which we also
denote by $\na$, by the identity (here $\w \in \W$, and $s\in \W^*(V,E)$)
\begin{equation}
\na( \w s)= \na(\w) s+  (-1)^{\deg \w} \w \na(s)
\end{equation}
One checks that the above formula indeed defines a degree 1~ derivation,
which can be described by the action on the
elements of the form $\al U_g$ where $\al \in \W^*(V, \End(E))$, $g
\in \G$, by the equation
\begin{equation}
\na (\al U_g)=\left(\na (\al)+\al \wedge \de(g)\right)U_g
\end{equation}
where $\de$ is $\W^1(V, \End(E))$-valued group cocycle, defined by
\begin{equation}
\de(g)=\na-g\circ \na \circ g^{-1}
\end{equation}
One defines a curvature as an element $\te U_1$, where 1~ is the
unit of the group, and $\te$ is the (usual) curvature of $\na$
. The graded trace $\dint$ on $\W^n$ is given by
\begin{equation}
\dint \al U_g =
 \begin{cases}
 \int \limits_V \al &\text{ if } g=1\\
 0                  &\text{ otherwise}
 \end{cases}
\end{equation}

One can associate with this cycle a cyclic  $n$-cocycle over an algebra
$\ca$, by the equation \eqref{ichern}.
By restricting it to the subalgebra $C_0^{\infty}(V)\rtimes \G$ one
obtains an $n$-cocycle $\{ \chi^k\}$ on this algebra. Its $k$-th component
is given by the formula
\begin{multline}
\chi^k(a_0U_{g_0}, a_1U_{g_1}, \dots a_kU_{g_k})=\\
\sum \limits_{1\leq i_1 < i_2  <\dots
< i_l \leq k} \int \limits_{V} a_0 da_1^{\g_1} da_2^{\g_2} \dots
da_{i_1-1}^{\g_{i_1-1}}
a_{i_1}^{\g_{i_1}}da_{i_1+1}^{\g_{i_1+1}}\dots\\
\Te_{i_1,i_2,\dots, i_l}(\g_1,\dots, \g_k)
\end{multline}
for $g_0g_1\dots g_k=1$ and 0~ otherwise.
Here the summation is over all the  subsets of $\{1,2,\dots
k\}$ and the following notations are used:
$\g_j$ are group elements defined by $\g_j=g_0g_1\dots g_{j-1}$.
$\Te_{i_1,i_2, \dots, i_l}(\g_1,\dots, \g_k)$ is the form (depending on $g_0,g_1$ \dots)
defined by the formula
\begin{multline}
\Te_{i_1,i_2, \dots, i_l}(\g_1,\dots, \g_k)=\\
\int \limits_{\De^k} tr
e^{-t_0 \te^{\g_1}} e^{-t_1 \te^{\g_2}}\dots e^{-t_{i_1-1}
\te^{\g_{i_1}}}\de(g_{i_1})^{\g_{i_1}}\\
e^{-t_{i_1} \te^{\g_{i_1+1}}} \dots e^{-t_{i_2-1}
\te^{\g_{i_2}}}\de(g_{i_2})^{\g_{i_2}}\dots  e^{-t_k \te}dt_1\dots dt_k
\end{multline}
The change of connection does not change the class in the cyclic
cohomology, as can be seen by constructing a a cobordism between
corresponding cycles.
This formula is a cyclic cohomological analogue of the formula of Bott \cite{bott78}.
More precisely, the following theorem holds:
\begin{thm} \label{phi}
Let $\Ch_{\G}(E) \in H^*(V \times_{\G} \rm{E}\G)$ be the
equivarint Chern character. Let $\Phi: H^*(V \times_{\G} \rm{E}\G)
\to HP^*(C_0^{\infty}(V)\rtimes \G)$ be the canonical imbedding,
constructed by Connes, cf. \cite{cn94}. Then
\[
\Phi\left(\Ch_{\G}(E)\right)=[\chi]
\]
\end{thm}
 Here $E$  pulls back to an equivariant bundle on $V \times \rm{E}\G$,
and then drops down to $V \times_{\G} \rm{E}\G$, and
the equivariant Chern character $\Ch_{\G}(E)$ is the Chern character
of the resulting bundle. We recall that we use normalization from \cite{bgv}. 

To prove the theorem we need some preliminary constructions and facts.
For a $\G$-manifold $Y$ by $Y_{\G}$ we denote the homotopy quotient
$Y \times_{\G} \rm{E}\G$.

Suppose we are given $\G$-manifolds $V$ and $X$, $X$ oriented.
 We construct then a map
$I: HC^j(C_0^{\infty}(V)\rtimes \G) \to
HC^{j+\dim X}(C_0^{\infty}(V\times X) \rtimes \G)$.
The construction is the following:
in $HC^{\dim X}(C_0^{\infty}(X)\rtimes \G)$ there is a class represented
by the cocycle
\[
\tau(f_0U_{g_0},f_1U_{g_1}, \dots, f_kU_{g_k})=
\begin{cases}
\int_X f_0 df_1^{g_0} \dots df_k^{g_0 g_1 \dots g_{k-1}} &\text{ if }
g_0g_1 \dots g_k =1\\
0 &\text { otherwise}
\end{cases}
\] 
One then constructs the map $I$ from the following diagram:
\begin{multline}
HC^j(C_0^{\infty}(V)\rtimes \G) \overset{ \cup \tau } \to
HC^{j +\dim X}\left(\bigl( C_0^{\infty}(V) \rtimes \G \bigr)\otimes \bigl(C_0^{\infty}(X) \rtimes
\G\bigr)\right)\\ \overset{\De^*}\to HC^{*+\dim X}(C_0^{\infty}(V\times X) \rtimes \G)
\end{multline}
Here the last arrow is induced by  the natural map
\begin{multline}
\De :
C_0^{\infty}(V\times X) \rtimes \G= \bigl( C_0^{\infty}(V)\otimes C_0^{\infty}(X)\bigr) \rtimes
\G
\\
\to
\bigl( C_0^{\infty}(V) \rtimes \G \bigr)\otimes \bigl(C_0^{\infty}(X) \rtimes
\G\bigr) 
 \end{multline}
defined by $\De \left((f\otimes f') U_g \right)= \left(fU_g \right)   \otimes \left(f'U_g \right)$.

Suppose now that $V$ is also oriented.
\begin{prop} \label{diag}
 The following diagram is commutative:
\begin{equation}
\begin{CD}
HP^*(C_0^{\infty}(V\times X) \rtimes \G) @<\Phi<<
 H^*\left( (V\times
X)_{\G} \right)\\
@A I AA  @AA \pi^* A\\
HP^*(C_0^{\infty}(V) \rtimes \G) @<\Phi<<
 H^*\left( V_{\G} \right)
\end{CD}
\end{equation}
Here $\pi : (V\times
X)_{\G} \to V_{\G}$ is induced by the ( $\G$-equivariant ) projection $V\times X \times \rm{E}\G
\to V\times \rm{E}\G$.
\end{prop}
\begin{proof}
 We can consider $V\times X$ with action of $\G \times \G$. We start
with showing that the following diagram is commutative:
\begin{equation}
\begin{CD}
HP^*(C_0^{\infty}(V\times X) \rtimes (\G\times \G)) @<\Phi<<
 H^*\left( (V\times
X)_{(\G\times \G)} \right)\\
@A \cup \tau AA  @AA \pi^* A\\
HP^*(C_0^{\infty}(V) \rtimes \G) @<\Phi<<
 H^*\left( V_{\G} \right)
\end{CD}
\end{equation}
Here we identify
 $C_0^{\infty}(V\times X) \rtimes (\G\times \G)$ with
$ (C_0^{\infty}(V) \rtimes \G)\otimes
(C_0^{\infty}(X) \rtimes \G)$ and $(V\times X)_{(\G\times \G)}$
with $X_{\G} \times V_{\G}$.
This is verified by the direct computation, using the Eilenberg-Silber
theorem and shuffle map in cyclic cohomology, cf. \cite{lo92}.

Now we note that the commutativity of the following diagram is clear:
\begin{equation}
\begin{CD}
HP^*(C_0^{\infty}(V\times X) \rtimes (\G\times \G)) @<\Phi<<
 H^*\left( (V\times
X)_{(\G\times \G)} \right)\\
@VVV  @VVV\\
HP^*(C_0^{\infty}(V\times X) \rtimes \G) @<\Phi<<
 H^*\left( (V\times
X)_{\G} \right)
\end{CD}
\end{equation}
where the both vertical arrows are induced by the diagonal maps $\G \to
\G \times \G$ and  $\rm{E}\G \to
\rm{E}\G \times \rm{E}\G$.
This ends the proof.
\end{proof}
\begin{prop} \label{funct}
 Let $E$ be an equivariant vector bundle on $V$ with
connection $\na$. Let $\chi \in HC^n(C_0^{\infty}(V) \rtimes \G)$ be the
character of the associated cycle, and let $\chi'\in
HC^{n+k}(C_0^{\infty}(V\times X) \rtimes \G)$ be the character of
the cycle constructed with the bundle $pr_V^*E$ and connection
$pr_V^*\na$, where $pr_V :X\times V \to V$.
Then $I(\chi)=\chi'$ (Here $n$ and $k$ are dimensions of $V$ and $X$ respectively)
\end{prop}
\begin{proof}
Let $\cc$ denote the corresponding cycle over $C_0^{\infty}(V) \rtimes
\G$, and $T$ -transverse fundamental cycle of $X$. Then $\cc \times T$
is a cycle over $(C_0^{\infty}(V) \rtimes \G)\otimes
(C_0^{\infty}(X) \rtimes \G)$, and
\[
\Ch(\cc \times T)=   \Ch(\cc) \cup \tau
\]
by the Corollary \ref{prod}.
If by $pr^* \cc$ we denote the corresponding cycle over
$C_0^{\infty}(V\times X) \rtimes \G$, we have
\[
\Ch(pr^* \cc)=\De^* \left(\Ch(\cc \otimes T) \right)=
\De^* \left(  \Ch(\cc) \cup \tau \right)= I( \Ch(\cc) )
\]
\end{proof}

\begin{lemma}\label{invariant}
 Suppose in addition to the conditions
of the Theorem \ref{phi} that $\G$ acts freely and properly on $V$.
 Then the statement of the Theorem holds.
 \end{lemma}
 \begin{proof}
Since the group acts freely and properly, one can find a connection on
$E$ which is $\G$-invariant. For the class of the cocycle $\chi$ written
with the invariant connection the result follows easily from the
definition of the map $\Phi$.
 \end{proof}

\begin{proof}[Proof of the Theorem \ref{phi}]
Comparison of the construction from \cite{ni90} with the definition of
the map $\Phi$ implies that (class of) $\chi$ is in the image
of $\Phi$, $[\chi]= \Phi(\xi)$ for some (necessarily unique)
$\xi \in H^*(V \times_{\G} \rm{E}\G)$. We need to verify that $\xi =\Ch_{\G}(E)$.
 We do this by showing that for any oriented
manifold $W$ and any map continuous $f:W \to V_{\G}$ $f^*\xi =f^* \Ch_{\G}(E)$.

Let $\tWW$ be the principal $\G$-bundle obtained by pullback of the
bundle $V \times \rm{E}\G \to V_{\G}$, so that the following diagram is
commutative, and $\tf$ is $\G$-equivariant:
\begin{equation}
\begin{CD}
\tWW @>\tf>> V \times \rm{E}\G\\
@VVV        @VVV\\
W @>f>>    V_{\G}
\end{CD}
\end{equation}
We can write $\tf$ as a composition of two $\G$-equivariant maps $\tf_1
:\tWW \to \tWW \times V \times \rm{E}\G$, which embeds $\tWW$ as the
graph of $\tf$ and $pr:\tWW \times V \times \rm{E}\G \to V \times
\rm{E}\G$, projection. Let $\pi: (\tWW \times V)_{\G} \to V_{\G}$ and
$f_1: W \to (\tWW \times V)_{\G}$ be the
induced maps. We have $f=\pi f_1$.

Construct now the class $\chi' \in HP^n(C_0^{\infty}(\tWW
\times V)\rtimes \G)$ using the bundle $pr^* E$ with connection $pr^*
\na$. By the Proposition \ref{funct} $\chi' =I(\chi)$, where $I: HP^*(C_0^{\infty}(V)\rtimes
\G) \to HP^{*+\dim W}(C_0^{\infty}(\tWW
\times V)\rtimes \G)$. By the Proposition \ref{diag} $\chi'=I(\chi)=I(\Phi(\xi))=\Phi(\pi^*\xi)$.
By the Lemma \ref{invariant}, since $\tWW\times V$ is acted by $\G$ freely and properly,
$\chi' = \Phi (\Ch(pr^*E))$. But since $\Ch(pr^*(E))=\pi^* \Ch(E)$, and
using injectivity of $\Phi$ we conclude that
\[
\pi^* \Ch(E)=\pi^*\xi
\]
Hence
\[
f^*\Ch(E)=f_1^*\pi^* \Ch(E)=f_1^*\pi^*\xi=f^*\xi
\]
  
\end{proof}

\begin{rem} [Relation with Connes' Godbillon-Vey cocycle]

In the paper \cite{cn86} Connes considers (in particular)
the case of the circle $S^1$ acted by the
 group of its diffeomorphisms $\Diff(S^1)$ .
Here we present the Connes construction in the multidimensional case and
indicate some relations with our construction of cyclic cocycles
representing  equivariant classes.

In the situation of the previous example take the bundle $E$ to
be $\bigwedge^n T^* X$. This is a 1~-dimensional trivial bundle,
naturally equipped with the action of the group $\G =\Diff (X)$. Let
$\phi$ be a nowhere 0~ section of this bundle, i.e. a volume form.
Define a flat connection $\na$ on $E$ by

\begin{equation}\label{con}
\na(f\phi)=df \phi, \ \f\in C^{\infty}(X)
\end{equation}
We can thus define the cycle $\cc$ over the algebra $C^{\infty}(X)$.

Let now $\de(g)$ be defined as above, and put
\begin{equation}
\mu(g)=\frac{\phi^g}{\phi} \in C^{\infty}(X)
\end{equation}
Then $\mu$ is a cocycle, i.e.
\begin{equation}
\mu(gh)=\mu(h)^g \mu(g)
\end{equation}

 We also have
\begin{equation}
\de(g)= d \log \mu(g)
\end{equation}
Indeed,
\[
\de(g) \phi^g = \na \phi^g - \left(\na (\phi)\right)^g= \na (\mu(g)
\phi)=d\mu(g) \phi
\]
and
\[
\de(g)=\frac{d \mu(g) \phi}{\phi^g}=\frac{d \mu(g) }{\mu(g)}= dlog(\mu(g))
\]

For every $t$ we define a homomorphism $\rho_t : C^{\infty}(X) \rtimes
\G \rightarrow \End(E)\rtimes \G$ by
\begin{equation}\label{flow}
\rho_t(aU_g) =a(\mu(g))^t U_g
\end{equation}
This is a homomorphism due to the cocycle property of $\mu$, which
according to \cite{cn86} is the Tomita-Takesaki flow associated with the
state given by the volume form $\phi$.

Consider now the transverse fundamental cycle $\Phi$ over the  algebra
$\ca = C^{\infty}(X)$
defined by the following data:

the differential graded algebra $\W^*(X) \rtimes \G$ with the
differential $d(\w U_g) =(d\w) U_g$

 the graded trace $\dint$
on $\W^*(X) \rtimes \G$ defined by
\[
\dint \w U_g =
 \begin{cases}
 \int \limits_X \w &\text{ if } g=1\\
 0                  &\text{ otherwise}
 \end{cases}
 \]

 the homomorphism $\rho=\rho_0=id$ from $\ca=C^{\infty}(X) \rtimes \G$ to
$C^{\infty}(X) \rtimes \G$.
The flow \eqref{flow} acts on the cycle $\Phi$, by replacing $\rho_0$ by
$\rho_t$. We call the cycle thus obtained $\Phi_t$. Using the identities
\begin{equation}
d(\rho_t(aU_g)) = \bigl(da +ta\ dlog \mu(g) \bigr)\mu(g)^t U_g
=\bigl(da +ta\ \de(g) \bigr)\mu(g)^t U_g
\end{equation}
and
\begin{equation}
\mu(g_0)\mu(g_1)^{g_0} \mu(g_2)^{g_0 g_1} \dots \mu(g_k)^{g_0 g_1 \dots
g_{k-1}}=\mu(g_0 g_1 \dots g_k) 
\end{equation}
we can  explicitly compute $\Ch(\Phi_t)$. This is the cyclic $n$-cocycle with
the only component of degree $n$
The result is:
\begin{equation} \label{cht}
\Ch(\Phi_t)=\sum_{j=0}^n t^j p_j
\end{equation}
where  $p_j$ is the cyclic cocycle given by
\begin{multline} \label{defp}
p_j(a_0U_{g_0},a_1U_{g_1}, \dots , a_nU_{g_n})=\\
\frac{1}{n!} \sum_{1\leq i_1 <i_2<\dots <i_j\leq n}
\int \limits_X
 a_0 da_1^{\g_1} da_2^{\g_2} \dots
da_{i_1-1}^{\g_{i_1-1}}
a_{i_1}^{\g_{i_1}}da_{i_1+1}^{\g_{i_1+1}}\dots\\
\Te_{i_1,i_2,\dots, i_j}(\g_1,\dots, \g_k)
\end{multline}
for $g_0g_1\dots g_k=1$ and 0~ otherwise, where we define as before
 $\g_j=g_0g_1\dots g_{j-1}$, and the $j$-form
$\Te_{i_1,i_2, \dots, i_j}(\g_1,\dots, \g_k)$ is given by
\begin{equation}
\Te_{i_1,i_2, \dots, i_j}(\g_1,\dots, \g_k)=
\de(g_{i_1})^{\g_{i_1}}\de(g_{i_2})^{\g_{i_2}}\dots  \de(g_{i_j})^{\g_{i_j}}
\end{equation}
In particular, $p_0$ is the transverse fundamental class.
Comparing these formulas from the formulas in the previous example we obtain
\begin{prop} Let $\Phi_1$ be the image of the transverse fundamental
cycle $\Phi$ under the action of the Tomita-Takesaki flow for the time 1~.
Let $\cc$ be the cycle over $C^{\infty}(X)\rtimes \G$ associated to the
 equivariant bundle  $\bigwedge^n T^* X$ with the connection from \eqref{con}.
 Then, on the level of cocycles
$\Ch(\Phi_1)=\Ch(\cc)$.
\end{prop}

We now sketch a construction
of a family of chains $\Psi_s$ providing the cobordism
between $\Phi_0$ and $\Phi_s$, $s \in \mathbb{R}$.
The algebra $\W^* =\W*([0,s])\gotimes \W^*(X)\rtimes \G$. The homomorphism
from $\ca$ to $\W^0$ maps $aU_g \in \W^*(X)\rtimes \G$ to
$a\mu(g)^tU_g$, where $t$ is the variable on $[0,s]$.
The connection is given by $1\gotimes \na +d \gotimes 1$ where $d$ is
the de Rham differential, and the curvature is 0~.
The restriction map is given by the restriction to the endpoints of
the interval and the graded trace is given by
\[
\dint \al \gotimes (\w U_g)= (-1)^{\deg\, \w} \int \limits_{[0,s]}
\al \int \limits_{X} \w
\]
if $\deg \, \al =1$ and $g=1$ and 0~ otherwise.
This chain  provides a cobordism between $\Phi_0$ and $\Phi_s$. Its
character is given by the formula
\begin{equation} \label{chp}
\Ch(\Psi_s)=\sum_{j=1}^{n+1} s^j q_j
\end{equation}
where  $q_j$ is the cyclic cochain given by
\begin{multline} \label{defq}
q_j(a_0U_{g_0},a_1U_{g_1}, \dots , a_nU_{g_n})=\\
\frac{1}{n!} \sum_{1\leq i_1 <i_2<\dots <i_j\leq n}
\int \limits_X
 a_0 da_1^{\g_1} da_2^{\g_2} \dots
da_{i_1-1}^{\g_{i_1-1}}
a_{i_1}^{\g_{i_1}}da_{i_1+1}^{\g_{i_1+1}}\dots\\
\Xi_{i_1,i_2,\dots, i_j}(\g_1,\dots, \g_k)
\end{multline}
for $g_0g_1\dots g_k=1$ and 0~ otherwise, where we define as before
 $\g_j=g_0g_1\dots g_{j-1}$, and the $j-1$-form
$\Xi_{i_1,i_2, \dots, i_j}(\g_1,\dots, \g_k)$ is given by
\begin{multline}
\Xi_{i_1,i_2, \dots, i_j}(\g_1,\dots, \g_k)=\\
\frac{1}{j}\sum \limits_{l=1}^{j}
(-1)^l \de(g_{i_1})^{\g_{i_1}}\de(g_{i_2})^{\g_{i_2}}\dots
\log \mu(g_{i_l})^{g_{i_l}} \dots \de(g_{i_j})^{\g_{i_j}}
\end{multline}
Comparing this formula with \eqref{cht} we obtain:
\begin{prop}
Let $p_j$, $j=1$,..., $n$ be the the chains, defined in \eqref{cht},
\eqref{defp}, 
and $q_j$, $j=1$,..., $n+1$ be from \eqref{chp}, \eqref{defq}. Then
for $j=1,\dots, n$ we have
\begin{equation}
B\,q_j=p_j \text{ and  } b\,q_j=0
\end{equation} 
Also
\begin{equation}
B\,q_{n+1}=0 \text{ and  } b\,q_{n+1}=0
\end{equation} 
In particular all $p_j$ define trivial classes in periodic
cyclic cohomology, and $q_{n+1}$ is a cyclic cocycle.
\end{prop}
The cocycle $q_{n+1}$ should represent (up to a constant) the Godbillon-Vey class
in the
cyclic cohomology (i.e. class defined by $h_1 c_1^n$, while $p_j$ and
$q_j$ represent forms $c_1^j$ and $h_1 c_1^j$, $j=1$,...,$n$, see
\cite{bott78}, but this remains to be verified.

\end{rem}

\begin{rem}[Transverse fundamental class]
The construction of the equivariant characteristic classes works equally
well in the case of a foliation. The new ingredient required here is the
Connes' construction of the transverse fundamental (generalized)
cycle. We now will
write a simple formula for the character of this cycle.

We start by briefly recalling Connes' construction from \cite{cn94}.
Details can be found in \cite{cn94}.
Let $(V,F)$ be a transversely oriented foliated manifold, $F$ being an integrable subbundle 
of $TV$. 
The graph of
the foliation $\cg$ is a groupoid, whose objects are points of $V$ and
morphisms are equivalence classes of paths in the leaves, with
equivalence given by holonomy. Equipped with a suitable topology it
becomes a smooth (possibly non-Hausdorff) manifold.
By $r$ and $s$ we denote the range and source maps $\cg \to V$.
By $\W_F^{1/2}$ we denote the line bundle on $V$ of the half-densities
in the direction of $F$.
  Let $\ca= C_0^{\infty}\left(\cg, s^*(\W_F^{1/2})\otimes
r^*(\W_F^{1/2})\right)$ be the convolution algebra
of $\cg$.
We define a ( nonunital ) generalized cycle over the algebra $\ca$ as
follows. The $k$-th component of the
graded algebra $\W^*$ is given by
$C_0^{\infty}\left(\cg, s^*(\W_F^{1/2})\otimes
r^*(\W_F^{1/2}) \otimes r^*(\bigwedge^k \tau^*)  \right) $. Here
$\tau = TV/F$ is the normal bundle, and the product $\W^k \otimes \W^l
\to \W^{k+l}$ is induced by the convolution and exterior product.

The definition of the transverse differentiation (connection)  requires a choice of a subbundle $H
\subset TV$, complementary to $F$. This choice allows one to identify
$C^{\infty}(V, \bigwedge ^*TV^*)$ with $C^{\infty}(V, \bigwedge
^*F^* \otimes \bigwedge^* \tau^*)$. We say that form $\w \in C^{\infty}(V, \bigwedge
^*TV^*)$ is of the type $(r,s)$ if it is in $C^{\infty}(V, \bigwedge
^r F^* \otimes \bigwedge^s \tau^*)$ under this identification.
For such a form we have
\begin{equation}
d\w = d_V \w +d_H \w +\s \w
\end{equation}
where $d_V \w$, $d_H \w$, $\s \w$ are defined to be components of $d \w$
of the types $(r+1, s)$, $(r, s+1)$, $(r-1, s+2)$ respectively (our notations are
slightly different from those of \cite{cn94}).
Now, writing locally $\rho \in  C^{\infty} (V, \W_F^{1/2})$ as $\rho = f |\w|^{1/2}$,
$f\in C^{\infty}(V)$, $\w \in C^{\infty}(V, \bigwedge ^{\dim F} F^*)$ we define
\begin{equation}
d_H \rho = (d_H f) |\w|^{1/2} + f |\w|^{1/2} \frac {d_H \w}{2 \w}
\end{equation}
Finally, $d_H$ can be extended uniquely as a graded derivation of the
graded algebra
 $C_0^{\infty}\left(\cg, s^*(\W_F^{1/2})\otimes
r^*(\W_F^{1/2}) \otimes r^*(\bigwedge^* \tau^*)  \right) $ so that the
following  identities are satisfied:
\begin{multline}
d_H \left(r^*(\rho_1) f s^*(\rho_2) \right)=\\
r^*(d_H\rho_1) f s^*(\rho_2) +
r^*(\rho_1) d_Hf s^*(\rho_2) +
r^*(\rho_1) f s^*(d_H\rho_2)\\
 \text{for } \rho_1, \rho_2 \in C^{\infty}(V,
\W_F^{1/2}), f \in C_0^{\infty}(\cg)
\end{multline}
and
\begin{multline}
d_H(\phi r^*(\w))=d_H(\phi) r^*(\w) +\phi r^*(d_H \w) \\
\text{for }\phi \in C_0^{\infty} \left(\cg, s^*(\W_F^{1/2}) \otimes
r^*(\W_F^{1/2}) \right), \w \in C^{\infty}(V,\bigwedge^*  \tau^*)
\end{multline}

Now, for the form $\w$ $d_H^2 \w =-(d_V\s+\s d_V)\w$. The operator
$\te= -(d_V\s+\s d_V)$ contains only longitudinal Lie derivatives, and
hence defines a multiplier (of degree 2) of the algebra
$C_0^{\infty}\left(\cg, s^*(\W_F^{1/2})\otimes
r^*(\W_F^{1/2}) \otimes r^*(\bigwedge^* \tau^*)  \right) $.

Finally, the graded trace on  $C_0^{\infty}\left(\cg, s^*(\W_F^{1/2})\otimes
r^*(\W_F^{1/2}) \otimes r^*(\bigwedge^q \tau^*)  \right)$, $ q=\codim F$ 
is given by $ \dint \w = \int \limits_V \w$.
\begin{lemma} (\cite{cn94}) $\left(C_0^{\infty}\left(\cg, s^*(\W_F^{1/2})\otimes
r^*(\W_F^{1/2}) \otimes r^*(\bigwedge^q \tau^*)  \right),
d_H, \te, \dint\right)$ is a generalized cycle of degree $q$
over the algebra $\ca$ .
\end{lemma}

We can now write an explicit formula for the character of this cycle. 
\begin{prop}
The following formula defines a (reduced)
cyclic cocycle $\chi$ in the $(b,$ $B)$-bicomplex
of the algebra $\ca$ (with adjoined unit).
\begin{equation}
\chi^k( \phi_0, \phi_1, \dots ,\phi_k)=
\frac{(-1)^{\frac{q-k}{2}}}{\left( \frac{q+k}{2} \right)!}
\sum \limits_{i_0+\dots +i_k= \frac{q-k}{2}}
\int \limits_V  \phi_0 \te^{i_0} d_H(\phi_1) \dots d_H (\phi_k) \te^{i_k}
\end{equation}
Here $k=q$, $q-2$, ..., and $\phi_j$, $j\geq 1$ are elements of $\ca$,
while $\phi_0$ is an element of $\ca$ with unit adjoined.

\end{prop}
Recall, that for $q$ even to define the cocycle over $\ca$ with the unit adjoined we
extend $\dint$ by requiring that $\dint \te ^{q/2}=0$.
The resulting class is independent of the choice of $H$.
It follows from the fact that by varying the subbundle $H$
smoothly we obtain the cobordism between the corresponding cycles,
satisfying the conditions of the Lemma \ref{varcon}. Note that the
equality here is in  cyclic cohomology, not only periodic cyclic cohomology. 

The results of the section \ref{cycles} imply that the class of
the cocycle $\chi$ is the transverse fundamental class of the foliation,
as defined in \cite{cn94}.  
\end{rem}

\section{Fredholm modules.} \label{appl}

In this section we
% associate the generalized cycle with a finitely
%summable bounded Fredholm module (cf. \cite{cn85}), and then
write
formulas for the character of the generalized cycle associated with a
finitely
summable bounded Fredholm module (cf. \cite{cn85}) .
In other words we obtain a formula for the character of a Fredholm module.
We show that this definition coincides with the Connes' definition \cite{cn85}.

 Let $(\ch,F,\g)$ be an even  finitely summable bounded
Fredholm module over the
algebra $\ca$. Here $\ch$ is a Hilbert space, on which the algebra $\ca$ acts,
 $\g$ is a
$\mathbb{Z}_2$-grading on $\ch$,  and $F$ is an odd selfadjoint operator
on $\ch$. We assume that $\ca$ is represented by the even operators in
$\ch$, and since we almost always consider only one representation of
$\ca$, we drop this representation from our notations, and do not
distinguish elements of the algebra and corresponding operators. We
suppose that the algebra $\ca$ is unital.
Let $p$ be a number such that $[F,a]\in \lp$ and $(F^2-1)\in
\mathcal{L}^{\frac{p}{2}}$. We remark that for any $p$ summable Fredholm 
   module one can
achieve these summability conditions by altering the operator $F$ and
keeping all the other data intact. We associate with the Fredholm module
a generalized cycle, similarly to \cite{cn85} where it is done in the
case when $F^2=1$.
  Consider a
$\mathbb{Z}$-graded algebra $\W=\bigoplus_{m=0}^{\infty}\W^m$ generated
by the symbols $a\in \ca$ of degree 0~, $[F,a]$, $a\in \ca$ of degree 1~
and symbol $(F^2-1)$ of degree 2~, with a relation $[F,ab]=a[F,b]+[F,a]b$
 This algebra can be naturally represented on
the Hilbert space $\ch$, and we will not distinguish
in our notations between elements of the algebra
and the corresponding operators.
 $\W$ is equipped with a natural
connection $\na$, given by the formula $\na(\xi)=[F,\xi]$
(graded commutator) in
terms of the  representation of $\W$, or on generators by the formulas
\begin{align}
&\na(a)=[F,a]\\
&\na([F,a])=\left((F^2-1)a-a(F^2-1)\right) =[(F^2-1),a]\\
&\na\left((F^2-1)\right)=0
\end{align}
Notice that $\na^2(\xi)=[(F^2-1),\xi]$ for $\xi \in \W$. Hence we 
define the curvature $\te$ to be $(F^2-1)$.
Clearly, $\xi \in \W^n$ is of trace class if $n\geq p$.
Here we  need to chose $n$ to be even, $n=2m$.
Hence we can
define the graded trace on $\W^n$ by $\dint \xi= m! Tr\,\g \xi$. The equality
 $Tr\,\g \na(\xi)=0$
for $\xi \in \W^{n-1}$  follows from the relation
\[
Tr\,\g \xi =\frac{1}{2}Tr\,\g F\na (\xi) -Tr\,\g (F^2-1)\xi
\]
which holds for $\xi$ of trace class). Indeed, for $\xi \in \W^{n-1}$
$\na(\xi)$ is of trace class and
\begin{multline}
Tr\,\g \na(\xi) =\frac{1}{2}Tr\,\g F\na^2 (\xi) +Tr\,\g \na(\xi)\\=
\frac{1}{2}Tr\,\g F[(F^2-1),\xi]-Tr\,\g (F^2-1)[F,\xi]=0
\end{multline}

Now we can apply the formula \eqref{dchern} to obtain a cyclic cocycle
$\Ch_{2m}(F)$ in the cyclic bicomplex of the algebra $\ca$. Its components
$\Ch^k_{2m}(F)$ $k=0$, 2~, 4~, \dots ,$2m$ are given by the formula
\begin{multline} \label {defc}
\Ch^k(F)(a_0,a_1,\dots a_k)=\\
\frac{m!}{(m+\frac{k}{2})!}\sum_{i_0+i_1+
\dots + i_k=m-\frac{k}{2}}Tr\,\g a_0 (1-F^2)^{i_0}[F, a_1](1-F^2)^{i_1}
\dots [F,a_k] (1-F^2)^{i_k}
\end{multline}
Note that for the case when $F^2=1$ we get the formula from \cite{cn85},
normalized as in \cite{cn94}.

We will now associate generalized chain with homotopy between Fredholm
modules. 
If the two Fredholm modules $(\ch,F_0,\g)$ and
$(\ch,F_1,\g)$   are connected by a smooth operator homotopy (
meaning that  there
exists a $C^1$ family $F_t$ of operators with $[F_t,a]\in \lp$ and $(F_t^2-1)\in
\mathcal{L}^{\frac{p}{2}}$, $t\in [0,1]$ with $F_t|_{t=0}=F_0$,
$F_t|_{t=1}=F_1$ ), this generalized chain will provide cobordism
between cycles corresponding to the modules.

We start by constructing, exactly as before, an algebra $\W_t$ generated by the
elements $a$, $[F_t,a]$, $(F_t^2-1)$, with the connection $\na_t$ and the
curvature $\te_t=(F_t^2-1)$. For each $t\in [0,1]$ one constructs a  natural
representation $\pi_t$ of
this algebra 
 on the Hilbert space $\ch$.
Let $\W^*([0,1])$ be the DGA  of
the differential forms on the interval $[0,1]$ with the usual
differential $d$. We can form a graded
tensor product $\W^*([0,1])\gotimes \W_t$.
Choose an odd number $n=2m+1$ so that $n \geq p+2$; if ain addition we
suppose that $\df \in \lp$, we can choose $n \geq p+1$.
In order to define the connection and the curvature we will have to adjoin to
our algebra an element  of degree 2~ $dt\gotimes \df$ and an element of
degree 3~ $dt\gotimes (F_t\df +\df F_t)$. The algebra with
the adjoined elements will be denoted $\W_c$.
The homomorphism $\rho_c :\ca \rightarrow \W_c$ is given by $\rho_c(a)=1\gotimes a$.
 We define the connection $\na_c$ as $
\frac{d}{dt} \wedge dt +\na_t$, i.e. on the generators the
definition is the following ($\bt \in \W^*([0,1])$ ):
\begin{align}
&\na_c(\bt \gotimes a)=d\bt \gotimes a+(-1)^{\deg(\bt)} \bt \gotimes [F_t,a]\\ 
&\na_c(\bt \gotimes [F_t,a])=d\bt \gotimes [F_t,a]+(-1)^{\deg(\bt)} \bt \gotimes [(F_t^2-1),a]+
\bt \wedge dt \gotimes [\df,a]\\
&\na_c(\bt \gotimes (F_t^2-1))=d\bt \gotimes (F_t^2-1)+\bt \wedge dt \gotimes
(F_t\df +\df F_t)\\
&\na_c(dt\gotimes \df)=- dt \gotimes (F_t\df +\df F_t)\\
&\na_c(dt \gotimes (F_t\df +\df F_t))=dt \gotimes [(F_t^2-1),\df]
\end{align}
The curvature $\te_c$ of this connection is defined as
\begin{equation}
\te_c=1\gotimes (F_t^2-1)+ dt\gotimes \df
\end{equation}
and the identity $(\na_c)^2 \cdot =[\te_c, \cdot]$ is verified by computation.
One
then defines 
the graded trace $\dint_c$ on $(\W^*([0,1])\gotimes \W_t)^n$
 by the formula
\[
\dint_c \bt \gotimes \xi =
 \begin{cases}
  (-1)^{\deg (\xi)} m! \int \limits_{[0,1]} \bigl(
  \bt \,Tr\, \g \pi_t(\xi) \bigr) &\text { if } \bt \in \W^1([0,1])\\
  0 &\text { if } \bt \in \W^0([0,1])
  \end{cases}
\]

The restriction maps $r_0 : \W_c \rightarrow \W_0$ and $r_1 : \W_c
\rightarrow \W_1$ are defined  as follows.  $r_0 (\bt \gotimes \xi))$ is 0~ if
$\bt$ is of degree 1~, and $\bt(0) \xi_0$ where $\xi_0$ is obtained from
$\xi$ by replacing $F_t$ by $F_0$ if $\bt$ is of degree 0~,
and similarly for $r_1$.
One can check that the map $r_1 \oplus r_0$ identifies  $\del \W_c$
with $\W^1 \oplus \widetilde{\W^0}$ and provides required cobordism.

Now we can use the Theorem \ref{bch} to study the properties of $\Ch(F)$
with respect to the operator homotopy.
\begin{thm} \label{homotop}
Suppose $(\ch,F_0,\g)$ and $(\ch,F_1,\g)$ are two finitely summable
Fredholm modules over an algebra $\ca$ which are connected by the smooth 
operator homotopy $F_t$ and $p$ is a number such that $[F_t]\in \lp$  and
$(F_t^2-1)\in \mathcal{L}^{\frac{p}{2}}$ for $0\leq t \leq 1$. Chose $m$
such that $2m\geq p+1$. Then $\Ch_{2m}(F_0)=\Ch_{2m}(F_1)$ in $HC^{2m}(\ca)$.
If moreover $\df \in \lp$  one can choose $m$ such
that $2m \geq p$.
\end{thm}
\begin{proof}
Let $Tch^k_{2m}$ denote the $k$-th component of the character of the
constructed above chain, providing the cobordism between the cycles
associated with $(\ch,F_0,\g)$ and $(\ch,F_1,\g)$, $k=1$, 3~,\dots,
$2m+1$.
It can be defined under the conditions on $m$ specified in the theorem.
According to the Theorem \ref{bch}
\[
\Ch_{2m}(F_1)-\Ch_{2m}(F_0)=(b+B)\, Tch_{2m}
\] 
 Now,
\[
Tch^{2m+1}_{2m}(a_0,a_1, \dots a_{2m+1}) =\text{const}\dint_c \rho_c(a_0) \na_c (\rho_c(a_1))
\dots \na_c (\rho_c(a_{2m+1}))=0
\]
(since the term under the $\dint_c $ does not contain $dt$). Hence
$Tch_{2m}$ can be considered as the $2m-1$ chain (is in the image of
$S$), and the result follows.
\end{proof}

\begin{rem} \label{ver}
Suppose we have two Fredholm modules $(\ch, F_0, \g)$ and $(\ch, F_1, \g)$
such that $F_0- F_1 \in \lp$ and $F_i^2-1 \in
\mathcal{L}^{\frac{p}{2}}$, $i=0$, 1~. Then $\Ch_{2m}(F_0)=\Ch_{2m}(F_1)$,
$ 2m \geq p$. Indeed, we can apply the Theorem\ref{homotop} to the
linear homotopy $F_t =F_0 +t (F_1-F_0)$, and need only to verify that
$F_t^2-1 \in \mathcal{L}^{\frac{p}{2}}$. But
\[
F_t^2-1 =(F_0^2-1)+ t\bigl(F_0(F_1-F_0)+(F_1-F_0)F_0 \bigr) +t^2(F_1-F_0)^2
\]
The first and the last terms in the right hand side are always in $
\mathcal{L}^{\frac{p}{2}}$, and since the left hand side is in
$\mathcal{L}^{\frac{p}{2}}$ for $t=1$, $\bigl(F_0(F_1-F_0)+(F_1-F_0)F_0
\bigr) \in \mathcal{L}^{\frac{p}{2}}$.
\end{rem}
\begin{cor}
Let $e$ be an idempotent in $M_N(\ca)$, and $(\ch,F,\g)$ be an even Fredholm 
module over $\ca$. Construct the Fredholm  operator 
$F_e = e(F \otimes 1) e : \ch^+\otimes \mathbb{C}^N \rightarrow \ch^- \otimes \mathbb{C}^N
$ (where $\ch^+$ and $\ch^-$ are determined by the grading).
Then 
\[
\ind(F_e)=<\Ch^*(F),\Ch_*(e)>
\]
Here $\Ch_*(e)$ is the usual Chern character in the cyclic homology.
\end{cor}
\begin{proof}
By replacing $\ca$ by $M_N(\ca)$ we reduce the situation to the case 
when $e \in \ca$. 
Now we apply  Connes' construction, which uses the
 homotopy $F_t=F+t(1-2e)[F,e]$ which connects
$F$ (obtained when $t=0$) with the operator $F_1=eFe +(1-e)F(1-e)$, obtained
when $t=1$. Note that $1- F_t^2 \in \mathcal{L}^{\frac{p}{2}}$. Indeed,
\[
F_t^2 -1=\bigl(F^2 -1\bigr)+\bigl(t(1-2e)[F,e] \bigr)^2 +t
\bigl([F,(1-2e)[F,e]] \bigr)
\]
The first two terms are clearly in $\mathcal{L}^{\frac{p}{2}}$. As for
the third one, it can be rewritten as $-2[F,e][F,e]+(1-2e)[F,[F,e]]=
-2[F,e]^2 +(1-2e)[(F^2-1),e] \in \mathcal{L}^{\frac{p}{2}}$.

The operator $F_1$ commutes with $e$, and homotopy does not change the
pairing. Hence it is enough to prove the result in the case when $F$ and $e$
commute.
In this case in the formula for the pairing
all of the terms involving commutators
are 0~, hence the only term with nonzero contribution is 
$\Ch^0(F)(e)=Tr\, \g e(1-F^2)^m = Tr\, \g (e-(eFe)^2)^m=\ind(F_e)$
by the well known formula.
\end{proof}

In \cite{cn85} Connes provides canonical construction, allowing one to
associate with every $p$-summable
Fredholm module such that $F^2-1 \neq 0$ another one for which $F^2-1=0$
, and which defines the same $K$-homology class. This allows to reduce
the definition of the character of a general Fredholm module to the case
when $F^2=1$.
The construction is the following. Given the Fredholm module $(\ch, F, \g)$ one
first constructs the Hilbert space $\tch = \ch \oplus \ch$ with the
grading given by $\tg=\g \oplus (-\g)$. An element $a \in \ca$ acts by
$
\left(
\begin{matrix}
a& 0\\
0& 0
\end{matrix}
\right)
$.
Then one constructs an operator $\tF$, such that $\tF -F' \in \lp$
and $\tF^2=1$; here by $F'$ we denote
$
\left(
\begin{matrix}
F& 0\\
0& -F\\
\end{matrix}
\right)
$.
The character of the Fredholm module $(\ch, F,\g)$ is then \emph{defined} to be
the character of the $(\tch, \tF, \tg)$.
\begin{thm}\label{coin}
Let $(\ch,F,\g)$ be an even  finitely summable Fredholm module over the
algebra $\ca$, and let $p$ be a real
number such that $[F,a]\in \lp$ and $(F^2-1)\in
\mathcal{L}^{\frac{p}{2}}$. Then
class of $\Ch^*(F)$ defined in \eqref{defc} in the periodic cyclic
cohomology 
 coincides with the  Chern character, as defined by Connes \cite{cn85}.
\end{thm}
\begin{proof}
First, let us consider the Fredholm module
$(\tch, F', \tg)$ over the algebra $\tca$ - the algebra $\ca$ with
adjoined unit ( acting by the identity operator). 
$\Ch_{2m}(F')$ is then defines a class in the cyclic cohomology of
$\tca$, where we choose $2m \geq p$.
Since $Tr\, \tg (1-(F')^2)^{m}=0$, it defines class in the reduced
cyclic cohomology of $\tca$, and hence in the cyclic cohomology of $\ca$.
It coincides with the class defined by the Fredholm module $(\ch, F, \g)$.

The Theorem \ref{homotop} and the Remark \ref{ver}
show that the classes defined by the
Fredholm modules $\Ch(\tF)$ and $\Ch(F')$ coincide.
To finish the proof we note that $\Ch(\tF)$ coincides with the Chern
character as defined in \cite{cn85}.

\end{proof}

The proof of the Theorem \ref{homotop} also provides an explicit transgression formula. We just
need to compute explicitly formula for
\begin{multline}
Tch^{k}_{2m}(a_0,a_1, \dots
a_{k})=\\
 \frac{(-1)^{m-\frac{k-1}{2}}(m)!}{(m+\frac{k+1}{2})!}\sum_{i_0+i_1+
\dots + i_k=m-\frac{k-1}{2}}\dint_c \rho_c(a_0) \te_c^{i_0} \na_c(
\rho_c(a_1))\te_c^{i_1}\dots \na_c(\rho_c(a_k)) \te_c^{i_k}
\end{multline}

Since $\te_c^{i_l}=\sum \limits_{r+q=i_l-1}dt \gotimes (F_t^2-1)^r \df (F_t^2-1)^q$ one can
rewrite this formula as
\begin{multline}
 \frac{(-1)^{m-\frac{k-1}{2}}(m)!}{(m+\frac{k+1}{2})!} \int_0^1 \Bigl( \sum_{i_0+i_1+
\dots + i_k=m-\frac{k-1}{2}} \sum_{l=0}^k \sum_{r+q=i_l-1} (-1)^l Tr\, \g
a_0 (F_t^2-1)^{i_0} \\
[F_t,a_1] (F_t^2-1)^{i_1}\dots [F_t,a_l](F_t^2-1)^r \df (F_t^2-1)^q\dots
[F_t,a_k](F_t^2-1)^{i_k} \Bigr)dt
\end{multline}
Finally we can write the answer as
\begin{multline}
Tch^{k}_{2m}(a_0,a_1, \dots
a_{k})=\\
- \frac{(m)!}{(m+\frac{k+1}{2})!} \int_0^1 \Bigl( \sum_{i_0+
\dots + i_k+i_{k+1}=m-\frac{k+1}{2}} \sum_{l=0}^k (-1)^l Tr\, \g
a_0 (F_t^2-1)^{i_0} \\
[F_t,a_1] (1-F_t^2)^{i_1}  \dots  [F_t,a_l](1-F_t^2)^{i_l} \df (1-F_t^2)^{i_{l+1}}\dots
[F_t,a_k](1-F_t^2)^{i_{k+1}} \Bigr)dt
\end{multline}
where $k$ is an odd number between 1~ and $2m-1$.

All the considerations above can be repeated in the case of an odd
finitely summable fredholm module $(\ch,F)$ over an algebra $\ca$ .
Here as before we  suppose that $[F,a]\in \lp$, $(F^2-1) \in
\mathcal{L}^{\frac{p}{2}}$. We choose number $m$ such that $n=2m+1 \geq p$.
The trace now is given by
$\dint \xi= \ci \G(n/2+1) Tr\, \xi$.

The corresponding Chern character $\Ch_{2m+1}(F)$ has components
$\Ch^k_{2m+1}$ for $k=1$, 3~, \dots, $2m+1$, given by the formula
\begin{multline}
\Ch^k_{2m+1}(a_0,a_1, \dots ,a_k)=
\frac{\G(m+\frac{3}{2})\ci}{(m+\frac{k+1}{2})!}\\
\sum_{i_0+i_1+
\dots + i_k=m-\frac{k-1}{2}}Tr\,  a_0(1-F^2)^{i_0} [F,a_1](1-F^2)^{i_1}
\dots [F,a_k](1-F^2)^{i_k}
\end{multline}

If the two Fredholm modules are connected via the operator homotopy
$F_t$ one has the transgression formula
\begin{equation}
\Ch_{2m+1}(F_1)-\Ch_{2m+1}(F_0)=(b+B)\, Tch_{2m+1}
\end{equation}
where $Tch_{2m+1}$ is a $2m$ cyclic cochain having components $Tch^{k}_{2m}$
for $k$ even between 0~ and $2m$, given by the formula:
\begin{multline}
Tch^{k}_{2m+1}(a_0,a_1, \dots
a_{k})=\\
- \frac{\G(m+\frac{3}{2})\ci}{(m+\frac{k}{2}+1)!} \int_0^1 \Bigl( \sum_{i_0+
\dots +i_k+i_{k+1}=m-\frac{k}{2}} \sum_{l=0}^k (-1)^l Tr\, 
a_0 (F_t^2-1)^{i_0} \\
[F_t,a_1] (1-F_t^2)^{i_1}  \dots  [F_t,a_l](1-F_t^2)^{i_l} \df (1-F_t^2)^{i_{l+1}}\dots
[F_t,a_k](1-F_t^2)^{i_{k+1}} \Bigr)dt
\end{multline}

The proof of the Theorem \ref{coin}
works in the odd situation as well
and shows that $\Ch^*(F)$ coincides with the  Chern character
as defined by Connes. In particular, this allows to recover the spectral
flow via the pairing with $K$-theory. More precisely, let $u \in
M_N(\ca)$ be a unitary. Let $\flow (F\otimes 1, (F\otimes 1)^u)$ be the
spectral flow of the operators $F\otimes 1$ and $(F\otimes 1)^u=u((F\otimes
1)u^*$ acting on the space
$\ch \otimes
\mathbb{C}^N $. The Chern character of the class of $u$ in $K_1(\ca)$ is
the periodic cyclic cycle defined by
 \begin{equation} \label{chu}
Ch_*(u)=\frac{1}{2\cp}\sum_{l=1}^{\infty}(-1)^l (l-1)! \,tr \left(
(u\otimes u^{-1})^l-( u^{-1}\otimes u)^l \right)
\end{equation}
Then we have the following
\begin{cor} \label{sf}
Let $u \in \ca$ be a unitary, and $(\ch, F)$ be an odd Fredholm
module over the algebra $\ca$. Then
$<\Ch^*(F), \Ch_*(u)>=\flow (F\otimes 1, (F\otimes 1)^u)$
\end{cor}
\begin{rem}
This is a finitely summable analogue of the result of Getzler \cite{ge93}
In the finitely summable case analytic formula for the spectral flow was
derived in \cite{cp98}; use of the Theorem \ref{coin} allows to give a proof
of the corollary \ref{sf}
without using this formula.
\end{rem}

%\bibliographystyle{plain}
%\bibliography{biblio}

\end{document}